\documentclass[11pt, oneside]{amsart}   	
\usepackage{geometry}                		
\usepackage{graphicx}				
\usepackage{amssymb}
\usepackage{amsmath}
\usepackage{amsthm}

\newcommand{\rea}{\mathbb{R}}

\newcommand{\BBZ}{\mathbb{Z}}

\newcommand{\calM}{\mathcal{M}}
\newcommand{\calQ}{\mathcal{Q}}
\newcommand{\calA}{\mathcal{A}}

\newcommand{\len}{\ell}
\newcommand{\dis}{\mathop{\mathrm{distance}}}

\def\<{\left\langle}
\def\>{\right\rangle}

\newtheorem{theorem}{Theorem}
\newtheorem{lemma}{Lemma}

\title{Sparse bounds for a prototypical singular Radon transform}
\author{Richard Oberlin}

\thanks{This material is based upon work supported by the National Science Foundation under Grant No. DMS-1440140 while the author was in residence at the Mathematical Sciences Research Institute in Berkeley, California, during the Spring 2017 semester.}

\begin{document}

\maketitle


\begin{abstract}
We use a variant of the technique in \cite{laceySB} to give sparse $L^p(\log(L))^4$ bounds for a class of model singular and maximal Radon transforms.
\end{abstract}

\section{Introduction} 

Suppose $\mu$ and $\sigma$ are finite signed and positive measures respectively,  supported on the unit ball $B(1) \subset \rea^n$ with 
$d\mu = \rho\ d\sigma$
for some bounded density $\rho$, $\mu(\rea^n) = 0$, and (using $\hat{\ }$ to denote the Fourier transform)
\begin{equation} \label{fsmoothing}
\max(|\hat{\sigma}(\xi)|,|\hat{\mu}(\xi)|) \lesssim |\xi|^{-\alpha}
\end{equation}
for some $\alpha > 0$ (Our main examples of interest are when $\sigma$ is surface measure on a compact piece of a finite-type submanifold of $\rea^n$ and $\rho$ is a smooth function on $\rea^n$ with $\sigma$-mean zero). Define $\mu_j$ by
\[
\int f \ d\mu_j = \int f(2^jx)\ d\mu(x).
\]
Given a collection of coefficients $\{\epsilon_j\}_{j \in \BBZ}$ with $|\epsilon_j| \leq 1$ we may consider the singular Radon transform
\[
T[f] := \sum_{j} \epsilon_j  \mu_j * f
\]
and the maximal averaging operator 
\[
T^*[f](x) := \sup_j \sigma_j * |f|(x).
\]
It is well known that condition \eqref{fsmoothing} implies that $T$ and $T^*$ are bounded on $L^p$ for $1 < p < \infty$.

The following ``sparse bound'' for $T^*$ was recently proven in \cite{laceySB} (see also related work \cite{cladekSB})
\begin{theorem}[Lacey] \label{laceythm}
Suppose $\sigma$ is surface measure on the unit sphere in $\rea^n$ and $1 < p < q < \infty$ are exponents such that convolution with $\sigma$ is a bounded operator from $L^{p}$ to $L^{q}.$ For $0 < \theta < 1$ let 
\[
\frac{1}{p_{\theta}} := \frac{1-\theta}{p} + \frac{\theta}{2} \text{\ \ and\ \ }\frac{1}{q_{\theta}} := \frac{1-\theta}{q} + \frac{\theta}{2}.
\]
There is a finite $C_{\theta}$ such that for every pair of compactly supported $f_1,f_2$ there is a sparse collection of cubes $\calQ$ such that 
\begin{equation} \label{asparsebound}
|\<T^*[f_1],f_2\>| \leq C_{\theta} \sum_{Q \in \calQ} |Q| \<|f_1|\>_{Q,p_{\theta}}\<|f_2|\>_{Q,q_\theta'}
\end{equation}
where
\[
\<|f|\>_{Q,p} := \left(\frac{1}{|Q|} \int_{Q} |f|^p\ dx \right)^{1/p}.
\]
\end{theorem}

Above, we use $|Q|$ to denote the Lebesgue measure of $Q$, and the collection $\calQ$ is said to be  \emph{sparse} if there 
is a collection of pairwise disjoint sets $\{F_Q\}_{Q \in \calQ}$ with  $|F_Q| \geq \frac{1}{2} |Q|$ and $F_Q \subset \calQ$. Bounds such as \eqref{asparsebound} (as well as those which give pointwise or norm domination by sparse operators) have been of much recent interest. See for example \cite{lernerAP}, \cite{lernerID}, \cite{lernerOP}, \cite{diplinioPS}, \cite{beneaSB}, \cite{bernicotSW}, \cite{culiucSB}, \cite{culiucDM}, \cite{laceyAE}, \cite{krauseSB}, \cite{nazarovCB}.

Theorem \ref{laceythm} is nontrivial (given that $T^*$ is known to be bounded on $L^p$) since $q_{\theta}' < p_{\theta}'.$ Furthermore, the range of exponents is sharp up to the small $\theta$-loss in interpolation (Since there is positive distance between the center of the sphere and the support of the measure, a sparse bound as above implies that convolution with $\sigma$ is bounded from $L^{p_{\theta}}$ to $L^{q_{\theta}}$). Lacey's argument does not appear to depend on the geometry of the sphere, and likely extends without modification to compactly supported positive measures satisfying \eqref{fsmoothing}. 

Our purpose here is to explore the relationship between the method of \cite{laceySB} and more traditional approaches (which use a regularization of the single scale operator) for bounding $T^*$. This will allow us to push a little closer to the natural endpoint exponents $(p,q).$ We have also organized our argument\footnote{Specifically, we use a Calder\'{o}n-Zygmund decomposition of \emph{both} functions, as was done in the original version of \cite{laceySB}. Later versions feature a streamlined argument which relies instead on the orthogonality of the linearizing functions and does not seem to immediately bound $T$.} to facilitate
 bounds for the singular integral $T$.
 
Given a cube $Q$, define 
\begin{equation} \label{plusdef1}
f^{0}_Q := f \cdot 1_{\{x \in Q : |f(x)| \leq \<|f|\>_{Q,p}\}} \text{\ \ and\ \ }f^{m}_Q := f \cdot 1_{\{x \in Q : 2^{m-1} \<|f|\>_{Q,p} <  |f(x)| \leq 2^m \<|f|\>_{Q,p}\}}, \ \  m > 0.
\end{equation}
Our bounds will be in terms of the following ``restricted-type $L^p\log(L)^4$" averages:
\begin{equation}
\<|f|\>_{Q,p^+} := \sum_{m \geq 0} (m+1)^4 \<|f^m_Q|\>_{Q,p}.
\end{equation}
It is straightforward to check that for each $\tilde{p} > p \geq 1$
\[
\<|f|\>_{Q,p}  \leq \<|f|\>_{Q,p^+} \leq C_{\tilde{p}} \<|f|\>_{Q,\tilde{p}}.
\]

\begin{theorem} \label{singulartheorem}
Suppose $\mu,\sigma$ are finite signed and positive measures respectively supported on the unit ball with $\mu(\rea^n) = 0.$ If  $\mu$ and $\sigma$ satisfy \eqref{fsmoothing} and $1 < p < q < \infty$ are exponents such that convolution with $\mu$ is a bounded operator from $L^{p}$ to $L^{q}$ then there is a finite $C$ such that for every pair of compactly supported functions $f_1,f_2$ there is a sparse collection of cubes $\calQ$ such that 
\begin{equation} \label{singularplusbound}
|\<T[f_1],f_2\>| \leq C \sum_{Q \in \calQ} |Q| \<|f_1|\>_{Q,p^+}\<|f_2|\>_{3Q,q'^+}.
\end{equation}
\end{theorem}

Essentially the same proof can be used to bound the maximal operator.

\begin{theorem} \label{maximaltheorem}
Suppose $\sigma$ is a finite measure supported on the unit ball satisfying \eqref{fsmoothing},  and that $1 < p < q < \infty$ are exponents such that convolution with $\sigma$ is a bounded operator from $L^{p}$ to $L^{q}.$ There is a finite $C$ such that for every pair of compactly supported $f_1,f_2$ there is a sparse collection of cubes $\calQ$ such that 
\begin{equation} \label{maximalplusbound}
|\<T^*[f_1],f_2\>| \leq C \sum_{Q \in \calQ} |Q| \<|f_1|\>_{Q,p^+}\<|f_2|\>_{3Q,q'^+}.
\end{equation}
\end{theorem}

The exponent four in the definition of $\<|f|\>_{Q,p^+}$ is not optimal and could be lowered slightly by following the numerology more carefully. We conjecture (based on parallels in the methods of proof) that the sharp bounds for \eqref{singularplusbound} and \eqref{maximalplusbound} may match the (currently unknown) sharp estimates at $L^1$ for $T$ and $T^*.$ Specifically, that for a given $\sigma,$  \eqref{maximalplusbound} should hold with  $\<|f_1|\>_{Q,p}\<|f_2|\>_{3Q,q}$ in place of $\<|f_1|\>_{Q,p^+}\<|f_2|\>_{3Q,q'^+}$ if and only if $T^*$ satisfies a weak-type $L^1$ estimate (and similarly for bounds with logarithmic losses). This would suggest that, at the very least, Theorems \ref{singulartheorem} and \ref{maximaltheorem} should hold with $L^p\log(L)$ in place of $L^p\log(L)^4.$

\subsection*{Acknowledgements} The author would like to thank Michael Lacey for sparking his interest in sparse bounds and for several helpful conversations subsequently.

\section{A review of the $L^p$ theory} \label{Lpsection}
We now quickly recall a standard method for proving $L^p$ estimates for $T$ (and $T^*$). This section is purely expository and may be skipped by the experts. 

The $L^2 \rightarrow L^2$ bound for $T$ is immediate from \eqref{fsmoothing}. To prove a bound near $L^1$, perform a Calder\'{o}n-Zygmund decomposition 
\[
f = g + \sum_{Q}b_Q
\]
where $b_Q$ is supported on the cube $Q$ and has mean-zero.
The contribution from the good function $g$ is handled, as usual, using the $L^2$ estimate. 

Let $\len(Q)$ denote the sidelength of a cube $Q$. If $\mu$ had an integrable derivative, we could deduce a weak-type $L^1$ estimate by leveraging the smoothness of the $\mu_j$ at scale $2^j$ against the cancellation of $b_Q$ for $2^j \geq \len(Q)$, and by using the decay of the $\mu_j$ at scale $2^j$ against the support of $b_Q$ for $2^j \leq \len(Q)$ (this, of course, is just the classic Calder\'{o}n-Zygmund method).

In general, one can write 
\[
\mu = \sum_{k \leq 0} \mu*\eta_k
\]
where $\mu*\eta_k$ is smooth at scale $2^k.$ Then $(\mu*\eta_k)_j$ is smooth at scale $2^{j+k}$, and so the contribution from $b_Q$ is acceptable, as above, when $2^{j+k} \geq \len(Q).$ Here, however, $(\mu*\eta_k)_j$ only has decay at scale $2^j$ and so, other than the trivial bound (i.e. the $(\mu*\eta_k)_j$ are uniformly in $L^1$ and so each of them gives a bounded convolution operator on $L^p$), one is not left with an obvious good option for $\len(Q) < 2^j < 2^{-k}\len(Q).$ This gives a weak-type estimate
\begin{equation} \label{L1boundblowup}
\|T^k[f]\|_{L^{1,\infty}} \lesssim (1-k) \|f\|_{L^1}
\end{equation}
where
\[
T^k[f] := \sum_j \epsilon_j(\mu*\eta_k)_j*f.
\]
On the other hand, provided $\eta_k$ is chosen with appropriate cancellation \eqref{fsmoothing} implies
\begin{equation} \label{L2bounddecay}
\|T^k[f]\|_{L^2} \lesssim 2^{k\alpha} \|f\|_{L^2}.
\end{equation}
Then $T$ is bounded on $L^p$ for $1 < p < \infty$ from the Marcinkiewicz interpolation theorem. It is not difficult, also using real interpolation, to do a little better (the following is only meant for illustration, and we omit its proof):

\begin{lemma}
Suppose $\{T_k\}_{k \leq 0}$ is \emph{any} sequence of operators satisfying \eqref{L1boundblowup} and \eqref{L2bounddecay}. Then\footnote{It is only coincidence that the four here matches the four in the definition of $\<|f|\>_{Q,p^+}$} for $r > 4$
\[
T := \sum_{k \leq 0} T^k
\]
satisfies the ``weak-type $L(\log(L))^r$ estimate"
\begin{equation} \label{weaktypeLlogL}
|\{|T[f]| > \lambda\}| \lesssim \int  \frac{|f(x)|}{\lambda} \left(\log\left(e + \frac{|f(x)|}{\lambda}\right)\right)^r\ dx, \ \ \ \lambda > 0.
\end{equation}
\end{lemma}

In fact, by incorporating the interpolation into the proof rather than crudely using it as a black-box, one finds that our operator $T$ satisfies a weak-type $L\log(L)$ bound, and for many measures $\mu$ one can apply more sophisticated techniques to push even closer to $L^1$.
 See, for example, \cite{seegerSM}, \cite{cladekIE}, and the references therein.

\section{Proof of Theorem \ref{singulartheorem}}
We will use a sparse bound adaptation (inspired by \cite{laceySB}) of the method outlined in Section \ref{Lpsection}. The principle use of the $L^p \rightarrow L^q$ estimate for convolution with $\mu$ is to replace the ``trivial $L^1$ bound" used for scales $\len(Q) < 2^j < 2^{-k}\len(Q)$ above.

Through a limiting argument and appropriate choice of dyadic grid, we may assume that there are finite $N_1,N_2$ such that $\epsilon_j = 0$ for $j$ outside of $[N_1,N_2]$ and that $f_1,f_2$ are supported on $Q_0$ and $3Q_0$ respectively, where $Q_0$ is a dyadic cube with $\len(Q_0) =2^{N_2} $(the bounds given will be independent of the $N_j$). Our proof will rely on recursion, each instance of which reduces $N_2$ and the support of the functions. After a finite number of steps, we are left with a null operator.

Write 
\begin{gather*}
\\ \calA_p^r[f](x) := \left(\frac{1}{|B(r)|} \int_{B(r)} |f(x + y)|^p\ dy\right)^{\frac{1}{p}}
\\ \calM_{p}[f](x) := \sup_{r > 0} \calA_p^r[f](x)
\\ T^*_{\text{high}}[f] :=  \sup_{2^j \leq \len(Q_0)} \sigma_j * |f|
\\ f_1^m := (f_1)^m_{Q_0} \ \ \ f_2^m := (f_2)^m_{3Q_0} \text{\ \ (using notation as in \eqref{plusdef1})}.
\end{gather*}
We then define
\begin{align*}
E_{1} = \{\calM_p&[f_1] > D \<|f_1|\>_{Q_0,p}\} 
\cup \{\calM_{1}[T^*_{\text{high}}[f_1]] > D \<|f_1|\>_{Q_0,p}\} 
\\ & \cup \bigcup_{m \geq 0} \{\calM_p[f_1^m] > (m+1) D \<|f_1^m|\>_{Q_0,p}\}.
\end{align*}
and similarly for $E_2$ with $f_2$  in place of $f_1$, $3Q_0$ in place of $Q_0,$ and $q'$ in place of $p$.

Choosing $D$ very large (depending on the $L^{p},L^{q'}$ bounds for $T^*$ and $\calM^p$), we can force $|E| := |E_{1} \cup E_{2}| \leq \frac{1}{2}|Q_0|$ and, say, $E \subset 6Q_0.$ 
Using a Whitney decomposition, write $E$ as the disjoint union of a  collection of dyadic cubes 
\[
E = \bigcup_{Q \in \calQ_{1}}Q
\]
each of which satisfies 
\begin{equation} \label{wreg}
5\sqrt{n}\len(Q) \leq \dis\left(Q,\left(\bigcup_{Q' \in \calQ_1}Q'\right)^c\right) < 11\sqrt{n}\len(Q).
\end{equation}
We then have, for example, that for every cube $Q'$ which contains a cube $Q \in \calQ_1$ 
\[
\<|f_1|\>_{Q',p} \lesssim \<|f_1|\>_{Q_0,p}.
\]

Perform a Calder\'{o}n-Zygmund decomposition of $f_1$
\begin{align*}
f_1 &=: g_1 + \sum_{Q \in \calQ_{1}}1_Q(f_1 - \<f_1\>_{Q,1}) =: g_1 + \sum_{Q \in \calQ_{1}}b_{1,Q}
\\ &=  g_1 + \sum_{\substack{Q \in \calQ_{1} \\ Q \subset Q_0}}b_{1,Q}
\end{align*}
where, for the last identity, we use that, since $Q_0 \not\subset E$, if $Q \cap Q_0 \neq \emptyset$ then $Q \subset Q_0.$
The good function is bounded
\[
\|g_1\|_{L^{\infty}} \lesssim \<|f_1|\>_{Q_0,p} .
\]
We will also use repeatedly that for any cube $Q'$ and $r \geq 1$
\[
\|\sum_{Q \subset Q'} b_{1,Q}\|_{L^r} \lesssim \|f_1\|_{L^r(Q')}.
\]

Decompose
\begin{equation} \label{firstgoodbad}
|\<T[f_1],f_2\>| \leq |\<T[g_1],f_2\>| + |\sum_{Q \in \calQ_{1}}\<T[b_{1,Q}],f_2\>|.
\end{equation}
The $L^q$ boundedness of $T$ implies that the first term on the right above
\[
|\<T[g_1],f_2\>| \lesssim  |Q_0| \<|f_1|\>_{Q_0,p}\<|f_2|\>_{3Q_0,q'}.
\]
Writing
\[
T_Q[f] := \sum_{2^j \leq \len(Q)} \epsilon_{j}\mu_j * (1_Q f),
\]
the second term of \eqref{firstgoodbad} 
\begin{equation} \label{firstbad}
 |\sum_{Q \in \calQ_{1}}\<T[b_{1,Q}],f_2\>| 
 = |\sum_{\substack{Q \in \calQ_{1} \\ Q \subset Q_0}}\<(T - T_{Q})[b_{1,Q}],f_2\> + \< T_{Q}[f_1],f_2\> - \<f_1\>_{Q,1}\<T_{Q}[1_{Q}],f_2\>|.
\end{equation}
By induction on $N_2 - N_1$, for each $Q \subset Q_0$ above we can find a sparse collection $\calQ_Q$ of dyadic subcubes of $Q$ such that
\[
|\< T_{Q}[f_1],f_2\>| \lesssim \sum_{Q' \in \calQ_Q} |Q'|\<|f_1|\>_{Q',p^+} \<|f_2|\>_{3Q',q'^+}.
\]
Setting $F_{Q_0} = Q_0 \setminus E$, we have that
\[
\calQ := \{Q_0\} \cup \bigcup_{\substack{Q \in \calQ_{1} \\ Q \subset Q_0}} \calQ_Q
\]
is sparse, and so it now remains to bound the sums of the first and third terms on the right of \eqref{firstbad}
\[
\lesssim  |Q_0| \<|f_1|\>_{Q_0,p^+}\<|f_2|\>_{3Q_0,q'^+}.
\]

Using the $L^q$ boundedness of $T_{Q}$  and the fact that the $3Q$ are finitely overlapping (from \eqref{wreg}), the sum of the third term is
\[
\lesssim \sum_{Q \in \calQ_{1}}  \<|f_1|\>_{Q_0,p} |Q|^{1/q} \|f_2\|_{L^{q'}(3Q)} \lesssim  |Q_0| \<|f_1|\>_{Q_0,p} \<|f_2|\>_{3Q_0,q'}. 
\]
The last, and main, step of the proof will be to show that
\begin{equation} \label{mainstep}
 |\sum_{Q \in \calQ_{1} }\<(T - T_{Q})[b_{1,Q}],f_2\>| \lesssim |Q_0| \<|f_1|\>_{Q_0,p^+} \<|f_2|\>_{3Q_0,q'^+}.
\end{equation} 
Perform a Calder\'{o}n-Zygmund decomposition of $f_2$
\[
f_2 =: g_2 + \sum_{Q \in \calQ_{1}}1_{Q}(f_2 - \<f_2\>_{Q,1}) =: g_2 + \sum_{Q \in \calQ_{1}}b_{2,Q}.
\]
The second good function is bounded
\[
\|g_2\|_{L^\infty} \lesssim \<|f_2|\>_{3Q_0,q'}
\]
which, using the $L^p$ boundedness of $T$  and $T_{Q}$ (separately), gives
\begin{align*}
  |\sum_{Q \in \calQ_{1}} \<(T  - T_{Q})[b_{1,Q}],g_2\>| 
  &\lesssim \|\sum_{Q \in \calQ_{1}} b_{1,Q} \|_{L^p}\|g_2\|_{L^{p'}} + \sum_{Q \in \calQ_{1}} \|b_{1,Q}\|_{L^p}\|g_2\|_{L^{p'}(3Q)}
 \\ &\lesssim |Q_0| \<|f_1|\>_{Q_0,p}\<|f_2|\>_{3Q_0,q'}.
\end{align*}
Expanding $T - T_Q$, \eqref{mainstep} will be finished once we estimate
\begin{equation} \label{mainstepbound}
| \sum_{Q,Q' \in \calQ_1} \sum_{2^j \geq \len(Q)} \epsilon_{j} \<\mu_j*b_{1,Q},b_{2,Q'}\>|.
\end{equation}
Then \eqref{mainstepbound} is 
\begin{multline} \label{squeezedj}
\leq |\sum_{Q,Q' \in \calQ_{1}}  \sum_{\substack{j \\ 2^j \geq \max(\len(Q),\len(Q'))}} \epsilon_j \<\mu_j*b_{1,Q},b_{2,Q'}\>| \\ + | \sum_{Q,Q' \in \calQ_{1}} \sum_{\substack{j \\ \len(Q)  \leq 2^j < \len(Q')}} \epsilon_j \<\mu_j*b_{1,Q},b_{2,Q'}\>|. 
\end{multline}
If a term in the right sum from \eqref{squeezedj} is nonzero then $Q \cap 2Q' \neq \emptyset$ and so, by \eqref{wreg}, $\len(Q) = \frac{1}{2}\len(Q') = 2^j.$ 
For each such $Q,Q'$, rescaling the $L^p \rightarrow L^q$ bound for $\mu$ gives
\[
|\<\mu_j*b_{1,Q},b_{2,Q'}\>| \lesssim |Q'| \<|f_1|\>_{3Q',p} \<|f_2|\>_{Q',q'} 
\]
and thus the right sum from \eqref{squeezedj} is 
\begin{align}
\nonumber&\lesssim \sum_{Q' \in \calQ_1} |Q'|\<|f_1|\>_{3Q',p} \<|f_2|\>_{Q',q'} 
\\ \nonumber&\leq\sup_{Q \in \calQ_1} \<|f_2|\>_{Q,q'}^{1 - \frac{q'}{p'}}  \sum_{Q' \in \calQ_1} |Q'|\<|f_1|\>_{3Q',p} \<|f_2|\>_{Q',q'}^{\frac{q'}{p'}} 
\\ \nonumber&\lesssim\sup_{Q \in \calQ_1} \<|f_2|\>_{Q,q'}^{1 - \frac{q'}{p'}}   |Q_0| \<|f_1|\>_{Q_0,p} \<|f_2|\>_{3Q_0,q'}^{\frac{q'}{p'}} 
\\ \label{interp} &\lesssim  |Q_0|\<|f_1|\>_{Q_0,p} \<|f_2|\>_{3Q_0,q'}.
\end{align}

We bound  the left sum from \eqref{squeezedj} by two terms which are treated in the same manner (it is irrelevant to the argument whether or not the diagonal $\len(Q) = \len(Q')$ is included), one of which is
\begin{equation} \label{halfsum}
|\sum_{\substack{Q, Q', j \\ \len(Q) \leq \len(Q') \leq 2^j}} \epsilon_j \<\mu_j*b_{1,Q},b_{2,Q'}\>|.
\end{equation}

It will be useful to decompose $\mu.$
Let  $\tilde{\eta}$ be a Schwartz function with $\hat{\tilde{\eta}}$ identically 1 on $B(1)$ and supported on $B(2)$
and $\eta := \tilde{\eta}_{-1} - \tilde{\eta}$ so that $\hat{\eta}$ is supported on $B(4) \setminus B(1)$ and 
\[
\hat{\tilde{\eta}} + \sum_{k \leq 0} \widehat{\eta_k} = 1.
\]
Then \eqref{halfsum}
\begin{align*}
 &\leq |\sum_{\substack{Q, Q', j \\ \len(Q) \leq \len(Q') \leq 2^j}} \epsilon_j \< (\tilde{\eta} *\mu)_j * b_{1,Q},b_{2,Q'}\>| 
\\& \ \ \ \ \ \ \ \ \ \ \ \ \ \ \ \ \ \ \ \ \ \ \ \ \  + \sum_{k \leq 0}|\sum_{\substack{Q, Q', j \\ \len(Q) \leq \len(Q') \leq 2^j}} \epsilon_j \< (\eta_k *\mu)_j * b_{1,Q},b_{2,Q'}\>| 
\\&=: |\tilde{S}| + \sum_{k \leq 0} |S_k|. 
\end{align*}

For $\tilde{S}$ we fix $Q$ and $2^j \geq \len(Q) =: 2^l$. Using the cancellation of $b_{1,Q}$ we have
\[
|\tilde{\eta}_j * b_{1,Q}(x)| \lesssim \<|f_1|\>_{Q_0,p} 2^{2(l - j)} (1 + |\dis(x,Q)|/2^j)^{-N}
\]
for large $N,$ giving (we will abuse notation by identifying $\mu$ with its conjugate reflection)
\begin{align*}
|\< (\tilde{\eta} *\mu)_j *  b_{1,Q},\sum_{\substack{Q' \\ \len(Q) \leq \len(Q') \leq 2^j}}b_{2,Q'}\> |
&\lesssim \<|f_1|\>_{Q_0,p} 2^{l - j} |Q| \calM_1[\mu_j*\sum_{\substack{Q' \\ \len(Q) \leq \len(Q') \leq 2^j}}b_{2,Q'}](x')
\\&\lesssim 2^{l - j} |Q| \<|f_1|\>_{Q_0,p} \<|f_2|\>_{3Q_0,q'}
\end{align*}
where $x' \in E^c$. (To obtain the second inequality above, we write $\sum b_{2,Q'}$ as the difference of $1_{\bigcup Q'} f_2$ and $\sum 1_{Q'}\<f_2\>_{Q',1}.$ The contribution from the former term is bounded by positivity of $\calM_1 \circ T^*$ and the fact that $x' \in E^c$, the contribution from the latter term instead uses the $L^\infty$ boundedness of $\calM_1[\mu_j * \cdot ].)$ Summing over $j$ and $Q'$ then gives 
\[
|\tilde{S}| \lesssim |Q_0| \<|f_1|\>_{Q_0,p} \<|f_2|\>_{3Q_0,q'}.
\]

We now fix $k \leq 0$ and turn our attention to $S_k$. We bound the  low frequency component
\[
\sum_{Q} \sum_{\substack{j \\ 2^j > 2^{-2k}\len(Q)}} |\sum_{\substack{Q' \\ \len(Q) \leq \len(Q') \leq 2^j}}\epsilon_j \< (\eta_k *\mu)_j * b_{1,Q},b_{2,Q'}\>|  
\lesssim 2^{k} |Q_0| \<|f_1|\>_{Q_0,p} \<|f_2|\>_{3Q_0,q'}
\]
using the same reasoning as for $\tilde{S}$ (and here, in contrast to $\tilde{S}$, it is important that $x' \in E^c$ since $u_j$ is at a coarser scale than $\eta_{k+j}$).

For $i=1,2$ write
\[
b^m_{i,Q} := 1_{Q}\left(f_i^m - \<f_i^m\>_{Q,1}\right).
\]
Since 
\[
f_i = \sum_{m \geq 0} f^m_i
\]
we have
\[
b_{i,Q} = \sum_{m \geq 0} b^m_{i,Q}.
\]
Decompose
\begin{align} \label{fixedkhf}
\sum_{Q} &\sum_{\substack{j \\ 2^j \leq 2^{-2k}\len(Q)}} \sum_{\substack{Q' \\ \len(Q) \leq \len(Q') \leq 2^j}}\epsilon_j \< (\eta_k *\mu)_j * b_{1,Q},b_{2,Q'}\> 
 \\ \nonumber &= \sum_{m_1,m_2 \geq 0} \sum_{Q} \sum_{\substack{j \\ 2^j \leq 2^{-2k}\len(Q)}} \sum_{\substack{Q' \\ \len(Q) \leq \len(Q') \leq 2^j}}\epsilon_j \< (\eta_k *\mu)_j * b^{m_1}_{1,Q},b^{m_2}_{2,Q'}\> 
\end{align}
For pairs $(m_1,m_2)$ with $m_1 + m_2 \leq \frac{-k \alpha}{2}$ we use the $L^2$ estimate for convolution with $
\eta_k*\mu.$ Writing
\[
Q_i^m := Q \cap \{f_i^m \neq 0\}
\]
for each $0 \leq h \leq \frac{-k\alpha}{2}$ and $0 \leq i \leq i' \leq -2k$ we have

\begin{align*} \label{fixedscalesum}
\sum_{m \leq h} \sum_l  &|\langle(\eta_k * \mu)_{l+i'}*\sum_{\len(Q) = 2^l} b^{m}_{1,Q},\sum_{\len(Q) = 2^{l+i}} b^{h-m}_{2,Q}\rangle|
\\ &\lesssim 2^{k\alpha }\sum_{m \leq h} \sum_l \|\sum_{\len(Q) = 2^l} b^{m}_{1,Q}\|_{L^2} \|\sum_{\len(Q) = 2^{l+i}} b^{h-m}_{2,Q}\|_{L^2}
\\ &\lesssim 2^{k\alpha }\sum_{m \leq h} \sum_l \|\sum_{\len(Q) = 2^l} 1_Qf_1^{m}\|_{L^2} \|\sum_{\len(Q) = 2^{l+i}} 1_Qf_2^{h-m}\|_{L^2}
\\ &\lesssim 2^{k\alpha + h} \<|f_1|\>_{Q_0,}  \<|f_2|\>_{3Q_0,q'}\sum_{m \leq h} \sum_l |\bigcup_{\len(Q) = 2^l}Q_{1}^m |^{\frac{1}{2}}  |\bigcup_{\len(Q) = 2^{l+i}}Q_{2}^{h-m}|^{\frac{1}{2}} 
\\ &\lesssim 2^{k \alpha + h}  |Q_0| \<|f_1|\>_{Q_0,p}  \<|f_2|\>_{3Q_0,q'}.
\end{align*}

Summing over $i,i'$ and then $h$ we have that the magnitude of the restriction of the sum on the right side of \eqref{fixedkhf} to $m_1 + m_2 \leq \frac{-k\alpha}{2}$ is 
\[
 \lesssim 2^{\frac{k \alpha}{4}}  |Q_0| \<|f_1|\>_{Q_0,p}  \<|f_2|\>_{3Q_0,q'} \]
which sums over $k \leq 0$ to an acceptable contribution.

For $m_1 + m_2 > \frac{-k\alpha}{2}$
we use the $L^p$ improving property of the $\mu$ averages. Fix $m_1,m_2$ and $0 \leq i \leq i' \leq 2k$. Then
\begin{align} \label{fixedscalesum2}
\sum_l  &|\langle(\eta_k * \mu)_{l+i'}*(\sum_{\len(Q) = 2^l} b^{m_1}_{1,Q}),(\sum_{\len(Q) = 2^{l+i}} b^{m_2}_{2,Q})\rangle|
\\ \nonumber &\lesssim \|f_2^{m_2}\|_{L^{q'}} (\sum_{l} \| \mu_{l+i'}*(\sum_{\len(Q) = 2^l}b^{m_1}_{1,Q})\|^q_{L^q} )^{1/q}.
\end{align}
The second factor on the right of \eqref{fixedscalesum2} is 
\begin{align*}
&\lesssim (\sum_{l} \sum_{\substack{Q ' \\ \len(Q') = 2^{l+i'} }} \| \mu_{l+i'}*(\sum_{\len(Q) = 2^l}b^{m_1}_{1,Q})\|^q_{L^q(Q')} )^{1/q}
\\ &\lesssim (\sum_{l} \sum_{\substack{Q ' \\ \len(Q') = 2^{l+i'} }} |Q'|\langle|\sum_{\len(Q) = 2^l}b^{m_1}_{1,Q}|\rangle_{3Q',p}^q )^{1/q}
\\ &\lesssim \sup_{\substack{Q'',l  \\ \len(Q'') = 2^{l+i'} }} \langle|\sum_{\len(Q) = 2^l}b^{m_1}_{1,Q}|\rangle_{3Q'',p}^{1 - \frac{p}{q}} \ \  (\sum_{l} \sum_{\substack{Q ' \\ \len(Q') = 2^{l+i'} }} |Q'|\langle|\sum_{\len(Q) = 2^l}b^{m_1}_{1,Q}|\rangle_{3Q',p}^p)^{1/q} 
\\ &\lesssim (m_1 + 1)\<|f_1^{m_1}|\>_{Q_0,p}^{1 - \frac{p}{q}} \|f_1^{m_1}\|_{L^p}^{\frac{p}{q}}
\\ &\lesssim (m_1 + 1)|Q_0|^{\frac{1}{q}} \<|f_1^{m_1}|\>_{Q_0,p}
\end{align*}
where, above, we sum over all dyadic cubes $Q'$ of sidelength $2^{l + i'}.$
This implies that the sum over $(i,i')$ of \eqref{fixedscalesum2} is 
\[
 \lesssim k^2 (m_1 + 1)|Q_0| \<|f_1^{m_1}|\>_{Q_0,p} \<|f_2^{m_2}|\>_{3Q_0,q'}
\]
and so the sum over $k$ of the magnitude of the restriction of the sum on the right side of \eqref{fixedkhf} to $m_1 + m_2 > \frac{-k\alpha}{2}$ is 
\begin{align*}
&\lesssim |Q_0|\sum_{m_1,m_2}(m_1 + m_2 + 1)^{4} \<|f_1^{m_1}|\>_{Q_0,p}  \<|f_2^{m_2}|\>_{3Q_0,q'}
\\ &\lesssim |Q_0|(\sum_{m}(m+1)^{4} \<|f_1^{m}|\>_{Q_0,p})( \sum_{m}(m+1)^{4} \<|f_2^{m}|\>_{3Q_0,q'})
\\ &= |Q_0|\<|f_1|\>_{Q_0,p^+} \<|f_2|\>_{3Q_0,q'^+}
\end{align*}
thus finishing the proof.

\bibliographystyle{halpha}  
\bibliography{ss}

\begin{thebibliography}{{Lac}17a}

\bibitem[BBL16]{beneaSB}
C.~{Benea}, F.~{Bernicot}, and T.~{Luque}.
\newblock {Sparse bilinear forms for Bochner Riesz multipliers and
  applications}.
\newblock {\em ArXiv e-prints}, May 2016, 1605.06401.

\bibitem[BFP16]{bernicotSW}
Fr\'ed\'eric Bernicot, Dorothee Frey, and Stefanie Petermichl.
\newblock Sharp weighted norm estimates beyond {C}alder\'on-{Z}ygmund theory.
\newblock {\em Anal. PDE}, 9(5):1079--1113, 2016.

\bibitem[CDO16]{culiucDM}
A.~{Culiuc}, F.~{Di Plinio}, and Y.~{Ou}.
\newblock {Domination of multilinear singular integrals by positive sparse
  forms}.
\newblock {\em ArXiv e-prints}, March 2016, 1603.05317.

\bibitem[CK17]{cladekIE}
L.~{Cladek} and B.~{Krause}.
\newblock {Improved endpoint bounds for the lacunary spherical maximal
  operator}.
\newblock {\em ArXiv e-prints}, March 2017, 1703.01508.

\bibitem[CKL16]{culiucSB}
A.~{Culiuc}, R.~{Kesler}, and M.~T. {Lacey}.
\newblock {Sparse Bounds for the Discrete Cubic Hilbert Transform}.
\newblock {\em ArXiv e-prints}, December 2016, 1612.08881.

\bibitem[CO]{cladekSB}
L.~{Cladek} and Y.~{Ou}.
\newblock {Sparse domination of Hilbert transforms along curves}.
\newblock {\em Forthcoming}.

\bibitem[DDU16]{diplinioPS}
F.~{Di Plinio}, Y.~Q. {Do}, and G.~N. {Uraltsev}.
\newblock {Positive sparse domination of variational Carleson operators}.
\newblock {\em ArXiv e-prints}, December 2016, 1612.03028.

\bibitem[KL17]{krauseSB}
B.~{Krause} and M.~T. {Lacey}.
\newblock {Sparse Bounds for Maximally Truncated Oscillatory Singular
  Integrals}.
\newblock {\em ArXiv e-prints}, January 2017, 1701.05249.

\bibitem[{Lac}17a]{laceySB}
M.~T. {Lacey}.
\newblock {Sparse Bounds for Spherical Maximal Functions}.
\newblock {\em ArXiv e-prints}, February 2017, 1702.08594.

\bibitem[Lac17b]{laceyAE}
Michael~T. Lacey.
\newblock An elementary proof of the {$A_2$} bound.
\newblock {\em Israel J. Math.}, 217(1):181--195, 2017.

\bibitem[Ler10]{lernerAP}
Andrei~K. Lerner.
\newblock A pointwise estimate for the local sharp maximal function with
  applications to singular integrals.
\newblock {\em Bull. Lond. Math. Soc.}, 42(5):843--856, 2010.

\bibitem[Ler16]{lernerOP}
Andrei~K. Lerner.
\newblock On pointwise estimates involving sparse operators.
\newblock {\em New York J. Math.}, 22:341--349, 2016.

\bibitem[LN15]{lernerID}
A.~K. {Lerner} and F.~{Nazarov}.
\newblock {Intuitive dyadic calculus: the basics}.
\newblock {\em ArXiv e-prints}, August 2015, 1508.05639.

\bibitem[NPTV17]{nazarovCB}
F.~{Nazarov}, S.~{Petermichl}, S.~{Treil}, and A.~{Volberg}.
\newblock {Convex body domination and weighted estimates with matrix weights}.
\newblock {\em ArXiv e-prints}, January 2017, 1701.01907.

\bibitem[STW04]{seegerSM}
Andreas Seeger, Terence Tao, and James Wright.
\newblock Singular maximal functions and {R}adon transforms near {$L^1$}.
\newblock {\em Amer. J. Math.}, 126(3):607--647, 2004.

\end{thebibliography}
\end{document}